\documentclass[12pt]{amsart}
\usepackage[cp850]{inputenc}
\usepackage{graphics}

\newtheorem{theorem}{Theorem}
\newtheorem{lemma}[theorem]{Lemma}

\newtheorem{corollary}[theorem]{Corollary}

\theoremstyle{definition}

\theoremstyle{remark}
\newtheorem{remark}[theorem]{Remark}

\usepackage{amssymb,amsmath}

\newcommand{\mbar}{\mbox{$\mathbb{M}^{3}_c$}}
\newcommand{\m}{\mbox{$\Sigma$}}
\newcommand{\mbarn}{\mbox{$\mathbb{M}^{n+1}_c$}}
\newcommand{\mn}{\mbox{$\Sigma^n$}}

\def\Rn{\mbox{$\mathbb R^{n+1}$}}
\def\Hn{\mbox{$\mathbb H^{n+1}$}}
\def\Sn{\mbox{$\mathbb S^{n+1}$}}
\newcommand{\Scal}{\mbox{$S$}}

\newcommand{\g}[2]{\mbox{$\langle #1 ,#2 \rangle$}}
\newcommand{\fle}{\mbox{$\rightarrow$}}
\newcommand{\rf}[1]{\mbox{(\ref{#1})}}
\newcommand{\rl}[1]{{~\ref{#1}}}

\newcommand{\xn}{\mbox{$\psi:\Sigma^n\fle\mbarn$}}

\def\beq{\begin{equation}}
\def\eeq{\end{equation}}

\begin{document}

\title[On the scalar curvature of CMC hypersurfaces]
{On the scalar curvature of constant mean curvature hypersurfaces in space forms}

\author{Luis J. Al\'\i as}
\address{Departamento de Matem\'{a}ticas, Universidad de Murcia, E-30100 Espinardo, Murcia, Spain}
\email{ljalias@um.es}
\thanks{This work was partially supported by MEC projects MTM2007-64504 and Fundaci\'{o}n S\'{e}neca
project 04540/GERM/06, Spain.
This research is a result of the activity developed within the framework of the Programme in Support of Excellence Groups of the Regi\'{o}n de Murcia, Spain, by Fundaci\'{o}n S\'{e}neca, Regional Agency
for Science and Technology (Regional Plan for Science and Technology 2007-2010).}

\author{S. Carolina Garc\'\i a-Mart\'\i nez}
\address{Departamento de Matem\'{a}ticas, Universidad de Murcia, E-30100 Espinardo, Murcia, Spain}
\email{sandracarolina.garcia@alu.um.es}
\thanks{S. Carolina Garc\'\i a-Mart\'\i nez was supported by a research training grant within the framework of the programme
Research Training in Excellence Groups GERM by Universidad de Murcia.}


\date{March 2009; revised version August 2009}



\begin{abstract}
In this paper we study the behavior of the scalar curvature \Scal\ of a complete hypersurface immersed with constant
mean curvature into a Riemannian space form of constant curvature, deriving a sharp estimate for
the infimum of \Scal. Our results will be an application of a weak Omori-Yau
maximum principle due to Pigola, Rigoli and Setti \cite{PRS}.
\end{abstract}

\maketitle

\section{Introduction}
\label{s1}
In a classical paper, Klotz and Osserman \cite{KO} characterized totally umbilical spheres and circular cylinders as the
only complete surfaces immersed into the Euclidean 3-space $\mathbb{R}^3$ with constant mean curvature $H\neq 0$ and whose
Gaussian curvature does not change sign. Later on, Hoffman \cite{Ho} and Tribuzy \cite{Tr} gave an extension of that
result to the case of surfaces with constant mean curvature in the Euclidean 3-sphere $\mathbb{S}^3$ and in the
hyperbolic space $\mathbb{H}^3$, respectively. Specifically, putting together the
results of those authors in a single statement, one gets the following result (see also \cite[Proposition 3.3]{Ch}).
\begin{theorem}
\label{th0}
Let $\m$ be a complete surface immersed into a 3-dimensional space form \mbar\ ($c=0,1,-1$) with constant mean curvature
$H$. If its Gaussian curvature $K$ does not change sign, then $\Sigma$ is
either a totally umbilical surface or $K=0$ and
\begin{enumerate}
\item[(a)] $c=0$ and \m\ is a circular cylinder
$\mathbb{R}\times\mathbb{S}^1(r)\subset\mathbb{R}^3$, with $r>0$,
\item[(b)] $c=1$ and \m\ is a flat torus
$\mathbb{S}^1(\sqrt{1-r^2})\times\mathbb{S}^1(r)\subset\mathbb{S}^3$, with $0<r<1$,
\item[(c)] $c=-1$ and \m\ is a hyperbolic cylinder
$\mathbb{H}^1(-\sqrt{1+r^2})\times\mathbb{S}^1(r)\subset\mathbb{H}^3$, with $r>0$.
\end{enumerate}
\end{theorem}
As a nice application of Theorem\rl{th0}, one gets the following consequence for the infimum of the Gaussian curvature
of $\Sigma$.
\begin{theorem}
\label{th1}
Let $\m$ be a complete surface immersed into a 3-dimensional space form \mbar\ ($c=0,1,-1$) with constant mean curvature
$H$ such that $H^2+c>0$, and let $K$ stand for its Gaussian curvature. Then
\begin{itemize}
\item[(i)] either $\inf_\Sigma K=H^2+c$, and \m\ is a totally umbilical surface,
\item[(ii)] or $\inf_\Sigma K\leq 0$, with equality if and only if
\begin{enumerate}
\item[(a)] $c=0$ and \m\ is a circular cylinder
$\mathbb{R}\times\mathbb{S}^1(r)\subset\mathbb{R}^3$, with $r>0$,
\item[(b)] $c=1$ and \m\ is a flat torus
$\mathbb{S}^1(\sqrt{1-r^2})\times\mathbb{S}^1(r)\subset\mathbb{S}^3$, with $0<r<1$,
\item[(c)] $c=-1$ and \m\ is a hyperbolic cylinder
$\mathbb{H}^1(-\sqrt{1+r^2})\times\mathbb{S}^1(r)\subset\mathbb{H}^3$, with $r>0$.
\end{enumerate}
\end{itemize}
\end{theorem}

Actually, it follows from the Gauss equation of the surface that $K\leq H^2+c$ on \m, with equality at the umbilical
points of $\Sigma$. Therefore, $\inf_\Sigma K\leq H^2+c$ with equality if and only if \m\ is totally umbilical. This
proves part (i). Moreover, if $\inf_\Sigma K<H^2+c$ then it must be $\inf_\Sigma K\leq 0$ necessarily. Otherwise, one
would have $K\geq \inf_\Sigma K>0$ which is not possible by Theorem\rl{th0}, since the non-totally umbilical surfaces in
(a), (b) and (c) are all flat. This shows that $\inf_\Sigma K\leq 0$. Finally, if equality holds, $\inf_\Sigma K=0$,
then $K\geq 0$ and the result follows from Theorem\rl{th0}.

Rotational surfaces show that the estimate in Theorem\rl{th1} is sharp. For instance, let us
consider the Delaunay rotational surfaces in the Euclidean space. For a given constant $H\neq 0$, we may consider the family of unduloids in $\mathbb{R}^3$
with constant mean curvature $H$, which are given by the following parametrization
\[
(s,\theta)\mapsto(x_B(s),y_B(s)\cos\theta,y_B(s)\sin\theta)
\]
where $0<B<1$ and
\begin{eqnarray*}
x_B(s) & = & \int_0^s\frac{1+B\sin{(2Ht)}}{\sqrt{1+B^2+2B\sin{(2Ht)}}}dt \\
{} & {} & {} \\
y_B(s) & = & \frac{\sqrt{1+B^2+2B\sin{(2Hs)}}}{2|H|}.
\end{eqnarray*}
The first fundamental form of these surfaces is $ds^2+y_B(s)^2d\theta^2$ and the Gaussian curvature is then
\[
K_B(s,\theta)=K_B(s)=-\frac{y''_B(s)}{y_B(s)}=\frac{4H^2B(B+\sin{(2Hs))(1+B\sin{(2Hs)})}}{(1+B^2+2B\sin(2Hs))^2}.
\]
Therefore, for these examples we have $\inf K_B=-{4H^2B}/{(1-B)^2}<0$, and for a given $\varepsilon>0$ there exists
$0<B<1$ such that $\inf K_B=-\varepsilon<0$.

It is worth pointing out that the proof of Theorem\rl{th0} (and hence Theorem\rl{th1}) strongly depends on the conformal
structure of the two-dimensional surface \m, and cannot be extended to higher dimensions. Our objective in this paper
is, using an alternative approach, to extend Theorem\rl{th1} to the case of $n$-dimensional hypersurfaces, with $n\geq 3$
(see Theorem\rl{th3aBIS} and Corollary\rl{coros} below).
As a consequence of Theorem\rl{th3aBIS}, in Corollary\rl{coro3a} we give a generalization of
Theorem 1.5 in \cite{AdC} to the case of complete parabolic hypersurfaces in space forms.

\subsection*{Acknowledgements} The authors would like to heartily thank A. Brasil Jr., E. de A. Costa and J.M. Espinar
for their valuable suggestions and useful comments during the preparation of this paper, which served to improve
the paper. They also thank the referee for valuable suggestions which improved the paper.

\section{Preliminaries}
\label{higherdimension}
Let us denote by \mbarn\ the standard model of an $(n+1)$-dimensional
Riemannian space form with constant curvature $c$, $c=0,1,-1$. That is, \mbarn\ denotes the Euclidean space \Rn\ when
$c=0$, the Euclidean sphere
\[
\Sn=\{x\in\mathbb{R}^{n+2} : |x|^2=1 \}\subset\mathbb{R}^{n+2},
\]
when $c=1$, and the hyperbolic space \Hn\ when $c=-1$. In this last case, it will be appropriate for us to use
the Minkowskian model of the hyperbolic space. Write $\mathbb{R}^{n+2}_1$ for $\mathbb{R}^{n+2}$, with canonical
coordinates $(x_0,x_1,\ldots,x_{n+1})$, endowed with the Lorentzian metric
\beq
\label{lor}
\g{}{}_1=-dx_0^2+dx_1^2+\cdots+dx_{n+1}^2.
\eeq
The $(n+1)$-dimensional hyperbolic space \Hn\ is the complete simply connected Riemannian
manifold with sectional curvature $-1$, which is realized as the
hyperboloid
\[
\Hn=\{ x\in{\mathbb R}^{n+2}_1 : \g{x}{x}_1=-1, x_{0}>0 \}\subset\mathbb{R}^{n+2}_1
\]
endowed with the Riemannian metric induced from ${\mathbb R}^{n+2}_1$.
In order to simplify our notation, when
$c=\pm 1$ we agree to denote by $\g{}{}$, without distinction, both the Euclidean metric on $\mathbb{R}^{n+2}$ and the
Lorentzian metric \rf{lor} on $\mathbb{R}^{n+2}_1$. We also agree to denote by $\g{}{}$ the corresponding Riemannian
metric induced on $\mbarn\hookrightarrow\mathbb{R}^{n+2}$.

Let us consider \xn\ an isometric immersion of an $n$-dimensional orientable Riemannian manifold $\Sigma$, and denote by
$A$ its second fundamental form (with respect to a globally defined normal unit vector field $N$) and
by $H$ its mean curvature, $H=(1/n)\mathrm{tr}(A)$. In the general $n$-dimensional case, instead of the curvature, it will be more appropriate to deal with the so called traceless second fundamental form of the hypersurface, which is given by $\Phi=A-HI$,
where $I$ denotes the identity operator on $\mathcal{X}(\Sigma)$. Observe that $\mathrm{tr}(\Phi)=0$ and $|\Phi|^2=\mathrm{tr}(\Phi^2)=|A|^2-nH^2\geq 0$, with equality if and only if $\Sigma$ is totally umbilical. For that reason, $\Phi$ is also
called the total umbilicity tensor of $\Sigma$.

As is well known, the curvature tensor $R$ of the hypersurface is given by the Gauss equation, which can be written in terms of $\Phi$ as
\begin{eqnarray}
\label{Gaussb}
\nonumber R(X,Y)Z & = & (c+H^2)(\g{X}{Z}Y-\g{Y}{Z}X)+\g{\Phi X}{Z}\Phi Y-\g{\Phi Y}{Z}\Phi X \\
{} & {} & +H(\g{\Phi X}{Z}Y-\g{Y}{Z}\Phi X+\g{X}{Z}\Phi Y-\g{\Phi Y}{Z}X)
\end{eqnarray}
for $X,Y,Z\in \mathcal{X}(\Sigma)$. In particular, the Ricci and the scalar curvatures of $\Sigma$ are given, respectively,
by
\beq
\label{Ricci}
\mathrm{Ric}(X,Y)=(n-1)(c+H^2)\g{X}{Y}+(n-2)H\g{\Phi X}{Y}-\g{\Phi X}{\Phi Y},
\eeq
for $X,Y\in \mathcal{X}(\Sigma)$, and
\beq
\label{scalar}
\Scal=n(n-1)(c+H^2)-|\Phi|^2.
\eeq

\subsection{Stochastic completeness and the Omori-Yau maximum principle}
For the proof of our results in higher dimension, we will make use of a weaker version of the Omori-Yau maximum principle.
Following the terminology introduced by Pigola, Rigoli and Setti in \cite{PRS}, the Omori-Yau maximum
principle is said to hold on an $n$-dimensional Riemannian manifold $\Sigma^n$ if, for any smooth function
$u\in\mathcal{C}^2(\Sigma)$ with $u^*=\sup_\Sigma u<+\infty$ there exists a sequence of points
$\{p_k\}_{k\in\mathbb{N}}$ in $\Sigma$ with the properties
\beq
\label{OY1}
(i) \quad u(p_k)>u^*-\frac{1}{k}, \quad (ii) \quad |\nabla u(p_k)|<\frac{1}{k}, \quad \mathrm{and} \quad (iii) \quad \Delta u(p_k)<\frac{1}{k}.
\eeq
In this sense, the classical result given by Omori and Yau in \cite{Om,Ya} states that the Omori-Yau maximum principle
holds on every complete Riemannian manifold with Ricci curvature bounded from below.
More generally, as shown by Pigola, Rigoli and Setti \cite[Example 1.13]{PRS}, a sufficiently controlled decay of the
radial Ricci curvature of the form
\beq
\label{radialRicci}
\mathrm{Ric}_\Sigma(\nabla\varrho,\nabla\varrho)\geq-C^2G(\varrho)
\eeq
where $\varrho$ is the distance function on \m\ to a fixed point, $C$ is a positive constant, and
$G:[0,+\infty)\rightarrow\mathbb{R}$ is a smooth function satisfying
\[
\textrm{(i)} \,\,\, G(0)>0, \,\,\, \textrm{(ii)} \,\,\, G'(t)\geq 0,  \,\,\, \textrm{(iii)} \,\,\,
\int_0^{+\infty}1/\sqrt{G(t)}=+\infty \textrm{ and }
\]
\[
\textrm{(iv)} \,\,\, \limsup_{t\rightarrow+\infty}tG(\sqrt{t})/G(t)<+\infty,
\]
suffices to imply the validity of the Omori-Yau maximum principle. In particular, and following the terminology
introduced by Bessa and Costa in \cite{BC}, the Omori-Yau maximum principle holds on a complete Riemannian
manifold whose Ricci curvature has \textit{strong quadratic decay} \cite{CX}, that is, with
\[
\mathrm{Ric}_\Sigma\geq-C^2(1+\varrho^2\log^2(\varrho+2)).
\]

On the other hand, as observed also by Pigola, Rigoli and Setti in \cite{PRS}, the validity
of Omori-Yau maximum principle on $\Sigma^n$ does not depend on curvature bounds  as much as one would expect. For
instance, the Omori-Yau maximum principle holds on every Riemannian manifold which is properly immersed into a
Riemannian space form with controlled mean curvature (see \cite[Example 1.14]{PRS}). In particular, it holds for
every constant mean curvature hypersurface properly immersed into a Riemannian space form.

More generally, and following again the terminology introduced in \cite{PRS}, the \textit{weak} Omori-Yau maximum
principle is said to hold on a (not necessarily complete) $n$-dimensional Riemannian manifold $\Sigma^n$ if, for any
smooth function  $u\in\mathcal{C}^2(\Sigma)$ with $u^*=\sup_\Sigma u<+\infty$ there exists a sequence of points
$\{p_k\}_{k\in\mathbb{N}}$ in $\Sigma$ with the properties
\beq
\label{OYw}
(i) \quad u(p_k)>u^*-\frac{1}{k},  \quad \mathrm{and} \quad (ii) \quad \Delta u(p_k)<\frac{1}{k}.
\eeq
As proved by Pigola, Rigoli and Setti \cite{PRS1}, the fact that the weak Omori-Yau maximum principle holds on $\Sigma$ is
equivalent to the stochastic completeness of the manifold
(see also \cite[Theorem 3.1]{PRS}).
In particular, the weak Omori-Yau maximum principle holds on every parabolic Riemannian manifold
(see also \cite[Corollary 6.4]{Gr}).

\section{Statement of the main results}
Now we are ready to state the following extension of Theorem\rl{th1} to the case of $n$-dimensional hypersurfaces, with
$n\geq 3$.
\begin{theorem}
\label{th3aBIS}
Let $\mn$ be a stochastically complete hypersurface immersed into an $(n+1)$-dimensional space form \mbarn\
($c=0,1,-1$ and $n\geq 3$) with constant mean curvature
$H$ such that $H^2+c>0$, and let \Scal\ stand for its scalar curvature. Then
\begin{itemize}
\item[(i)] either $\inf_\Sigma \Scal=n(n-1)(c+H^2)$ and \m\ is a totally umbilical hypersurface,
\item[(ii)] or
\[
\inf_\Sigma \Scal\leq\frac{n(n-2)}{2(n-1)}\left(2(n-1)c+nH^2+|H|\sqrt{n^2H^2+4(n-1)c}\right).
\]
Moreover, the equality holds and this infimum is attained at some point of $\Sigma$ if and only if
\begin{enumerate}
\item[(a)] $c=0$ and \m\ is an open piece of a circular cylinder
$\mathbb{R}\times\mathbb{S}^{n-1}(r)\subset\mathbb{R}^{n+1}$, with $r>0$,
\item[(b)] $c=1$ and \m\ is an open piece of either a minimal Clifford torus
$\mathbb{S}^k(\sqrt{k/n})\times\mathbb{S}^{n-k}(\sqrt{(n-k)/n})\subset\mathbb{S}^{n+1}$, with $k=1,\ldots,n-1$,
or a constant mean curvature torus
$\mathbb{S}^{1}(\sqrt{1-r^2})\times\mathbb{S}^{n-1}(r)\subset\mathbb{S}^{n+1}$, with
$0<r<\sqrt{(n-1)/n}$,
\item[(c)] $c=-1$ and \m\ is an open piece of a hyperbolic cylinder
$\mathbb{H}^{1}(-\sqrt{1+r^2})\times\mathbb{S}^{n-1}(r)\subset\mathbb{H}^{n+1}$, with $r>0$.
\end{enumerate}
\end{itemize}
\end{theorem}
In the particular case where $\Sigma^n$ is complete (which happens, for instance, when $\Sigma^n$ is properly immersed),
we obtain the following consequence.
\begin{corollary}
\label{coros}
Let $\mn$ be a complete hypersurface immersed into an $(n+1)$-dimensional space form \mbarn\
($c=0,1,-1$ and $n\geq 3$) with constant mean curvature $H$ such that $H^2+c>0$. Then
\begin{itemize}
\item[(i)] either $\inf_\Sigma \Scal=n(n-1)(c+H^2)$ and \m\ is a totally umbilical hypersurface,
\item[(ii)] or
\[
\inf_\Sigma \Scal\leq\frac{n(n-2)}{2(n-1)}\left(2(n-1)c+nH^2+|H|\sqrt{n^2H^2+4(n-1)c}\right).
\]
Moreover, the equality holds and this infimum is attained at some point of $\Sigma$ if and only if
\begin{enumerate}
\item[(a)] $c=0$ and \m\ is a circular cylinder
$\mathbb{R}\times\mathbb{S}^{n-1}(r)\subset\mathbb{R}^{n+1}$, with $r>0$,
\item[(b)] $c=1$ and \m\ is either a minimal Clifford torus
$\mathbb{S}^k(\sqrt{k/n})\times\mathbb{S}^{n-k}(\sqrt{(n-k)/n})\subset\mathbb{S}^{n+1}$, with $k=1,\ldots,n-1$,
or a constant mean curvature torus
$\mathbb{S}^{1}(\sqrt{1-r^2})\times\mathbb{S}^{n-1}(r)\subset\mathbb{S}^{n+1}$, with
$0<r<\sqrt{(n-1)/n}$,
\item[(c)] $c=-1$ and \m\ is a hyperbolic cylinder
$\mathbb{H}^{1}(-\sqrt{1+r^2})\times\mathbb{S}^{n-1}(r)\subset\mathbb{H}^{n+1}$, with $r>0$.
\end{enumerate}
\end{itemize}
\end{corollary}

On the other hand, it follows from \rf{scalar} that $\inf_\Sigma \Scal=n(n-1)(c+H^2)-\sup_\Sigma|\Phi|^2$. Therefore,
Theorem\rl{th3aBIS} (as well as Corollary\rl{coros})
can be re-written equivalently in terms of the total umbilicity tensor as follows.
\begin{theorem}
\label{th3a}
Let $\mn$ be a stochastically complete hypersurface immersed into an $(n+1)$-dimensional space form \mbarn\
($c=0,1,-1$ and $n\geq 3$) with constant mean curvature
$H$ such that $H^2+c>0$, and let $\Phi$ stand for its total umbilicity tensor.
Then
\begin{itemize}
\item[(i)] either $\sup_\Sigma|\Phi|=0$ and \m\ is a totally umbilical hypersurface,
\item[(ii)] or
\[
\sup_\Sigma|\Phi|\geq\alpha_H=\frac{\sqrt{n}}{2\sqrt{n-1}}\left(\sqrt{n^2H^2+4(n-1)c}-(n-2)|H|\right)>0.
\]
Moreover, the equality holds and this supremum is attained at some point of $\Sigma$ if and only if
\begin{enumerate}
\item[(a)] $c=0$ and \m\ is an open piece of a circular cylinder
$\mathbb{R}\times\mathbb{S}^{n-1}(r)\subset\mathbb{R}^{n+1}$, with $r>0$,
\item[(b)] $c=1$ and \m\ is an open piece of either a minimal Clifford torus
$\mathbb{S}^k(\sqrt{k/n})\times\mathbb{S}^{n-k}(\sqrt{(n-k)/n})\subset\mathbb{S}^{n+1}$, with $k=1,\ldots,n-1$,
or a constant mean curvature torus
$\mathbb{S}^{1}(\sqrt{1-r^2})\times\mathbb{S}^{n-1}(r)\subset\mathbb{S}^{n+1}$, with
$0<r<\sqrt{(n-1)/n}$,
\item[(c)] $c=-1$ and \m\ is an open piece of a hyperbolic cylinder
$\mathbb{H}^{1}(-\sqrt{1+r^2})\times\mathbb{S}^{n-1}(r)\subset\mathbb{H}^{n+1}$, with $r>0$.
\end{enumerate}
\end{itemize}
\end{theorem}

In particular, we get the following consequence, which gives a generalization of
Theorem 1.5 in \cite{AdC} to complete parabolic hypersurfaces in space forms.
\begin{corollary}
\label{coro3a}
Let $\mn$ be a complete parabolic hypersurface immersed into an $(n+1)$-dimensional space form \mbarn\
($c=0,1,-1$ and $n\geq 3$) with constant mean curvature
$H$ such that $H^2+c>0$, and let $\Phi$ stand for its total umbilicity tensor.
Then
\begin{itemize}
\item[(i)] either $\sup_\Sigma|\Phi|=0$ and \m\ is a totally umbilical hypersurface,
\item[(ii)] or
\[
\sup_\Sigma|\Phi|\geq\alpha_H=\frac{\sqrt{n}}{2\sqrt{n-1}}\left(\sqrt{n^2H^2+4(n-1)c}-(n-2)|H|\right)>0
\]
with equality if and only if
\begin{enumerate}
\item[(a)] $c=0$ and \m\ is a circular cylinder
$\mathbb{R}\times\mathbb{S}^{n-1}(r)$, with $r>0$,
\item[(b)] $c=1$ and \m\ is either a minimal Clifford torus
$\mathbb{S}^k(\sqrt{k/n})\times\mathbb{S}^{n-k}(\sqrt{(n-k)/n})$, with $k=1,\ldots,n-1$,
or a constant mean curvature torus
$\mathbb{S}^{1}(\sqrt{1-r^2})\times\mathbb{S}^{n-1}(r)$, with
$0<r<\sqrt{(n-1)/n}$,
\item[(c)] $c=-1$ and \m\ is a hyperbolic cylinder
$\mathbb{H}^{1}(-\sqrt{1+r^2})\times\mathbb{S}^{n-1}(r)$, with $r>0$.
\end{enumerate}
\end{itemize}
\end{corollary}

\section{Proof of the main results}
The proof of our results is based on a Simons type formula for the Laplacian of the function $|\Phi|^2$, which has already been used by several authors. For the sake of completeness, we include
here its derivation, following Nomizu and Smyth \cite{NS}. A standard tensor computation implies that
\beq
\label{eq2.1}
\frac{1}{2}\Delta|\Phi|^2=\frac{1}{2}\Delta\g{\Phi}{\Phi}=|\nabla \Phi|^2+\g{\Phi}{\Delta \Phi}.
\eeq
Here $\nabla\Phi:\mathcal{X}(\Sigma)\times \mathcal{X}(\Sigma)\fle \mathcal{X}(\Sigma)$ denotes the covariant differential of $\Phi$,
\[
\nabla\Phi(X,Y)=(\nabla_Y\Phi)X=\nabla_Y(\Phi X)-\Phi(\nabla_YX), \quad X,Y\in \mathcal{X}(\Sigma),
\]
and $\Delta\Phi:\mathcal{X}(\Sigma)\fle\ \mathcal{X}(\Sigma)$ is the rough Laplacian,
\[
\Delta\Phi(X)=\mathrm{tr}(\nabla^2\Phi(X,\cdot,\cdot))=\sum_{i=1}^n\nabla^2\Phi(X,E_i,E_i),
\]
where $\{E_1,\ldots, E_n\}$ is a local orthonormal frame on $\Sigma$.
Observe that, in our notation, $\nabla^2\Phi(X,Y,Z)=(\nabla_Z\nabla\Phi)(X,Y)$. Let us assume that the mean curvature $H$ is constant. In that case, $\nabla\Phi=\nabla A$, which is symmetric by the Codazzi equation of the
hypersurface and, hence,
$\nabla^2\Phi$ is also symmetric in its two first variables,
\[
\nabla^2\Phi(X,Y,Z)=\nabla^2\Phi(Y,X,Z), \quad X,Y,Z\in \mathcal{X}(\Sigma).
\]
With respect to the symmetries of $\nabla^2\Phi$ in the other variables, it is not difficult to see that
\[
\nabla^2\Phi(X,Y,Z)=\nabla^2\Phi(X,Z,Y)-R(Z,Y)\Phi X+\Phi(R(Z,Y)X).
\]
Thus, using the Gauss equation \rf{Gaussb} it follows from here that
\begin{eqnarray}
\label{laplacianPhi}
\Delta\Phi(X) & = & \sum_{i=1}^n\left(\nabla^2\Phi(E_i,E_i,X)-R(E_i,X)\Phi E_i+\Phi(R(E_i,X)E_i)\right) \\
\nonumber {} & = & \mathrm{tr}(\nabla_X(\nabla\Phi))-H|\Phi|^2X+(n(c+H^2)-|\Phi|^2)\Phi X+nH\Phi^2X \\
\nonumber {} & = & -H|\Phi|^2X+(n(c+H^2)-|\Phi|^2)\Phi X+nH\Phi^2X,
\end{eqnarray}
where we have used the facts that trace commutes with $\nabla_X$ and that $\mathrm{tr}(\nabla\Phi)=0$.
Therefore, by \rf{eq2.1} we conclude that
\beq
\label{Simons}
\frac{1}{2}\Delta|\Phi|^2=|\nabla \Phi|^2+nH\mathrm{tr}(\Phi^3)-|\Phi|^2(|\Phi|^2-n(c+H^2)).
\eeq
We will also need the following auxiliary result, known as Okumura lemma, which can be found in \cite{Ok} and
\cite[Lemma 2.6]{AdC}.
\begin{lemma}
\label{okumura}
Let $a_1,\ldots,a_n$ be real numbers such that $\sum_{i=1}^na_i=0$. Then
\[
-\frac{n-2}{\sqrt{n(n-1)}}(\sum_{i=1}^na_i^2)^{3/2}\leq
\sum_{i=1}^na_i^3\leq\frac{n-2}{\sqrt{n(n-1)}}(\sum_{i=1}^na_i^2)^{3/2}.
\]
Moreover, equality holds in the right-hand (respectively, left-hand) side if and only if $(n-1)$ of the $a_i$'s are
nonpositive (respectively, nonnegative) and equal.
\end{lemma}

\subsection{Proof of Theorem\rl{th3a} (or, equivalently, Theorem\rl{th3aBIS})}
Since $\mathrm{tr}(\Phi)=0$, we may use Lemma\rl{okumura} to estimate
$\mathrm{tr}(\Phi^3)$ as follows
\[
|\mathrm{tr}(\Phi^3)|\leq\frac{n-2}{\sqrt{n(n-1)}}|\Phi|^3,
\]
and then
\[
nH\mathrm{tr}(\Phi^3)\geq-n|H||\mathrm{tr}(\Phi^3)|\geq-\frac{n(n-2)}{\sqrt{n(n-1)}}|H||\Phi|^3.
\]
Using this in \rf{Simons}, we find
\begin{eqnarray}
\label{EQ7}
\frac{1}{2}\Delta|\Phi|^2 & \geq & |\nabla\Phi|^2-\frac{n(n-2)}{\sqrt{n(n-1)}}|H||\Phi|^3-|\Phi|^2(|\Phi|^2-n(c+H^2))\\
\nonumber {} & \geq & -|\Phi|^2P_H(|\Phi|),
\end{eqnarray}
where
\[
P_H(x)=x^2+\frac{n(n-2)}{\sqrt{n(n-1)}}|H|x-n(c+H^2).
\]
Observe that, since $H^2+c>0$, the polynomial $P_H(x)$ has a unique positive root given by
\[
\alpha_H=\frac{\sqrt{n}}{2\sqrt{n-1}}\left(\sqrt{n^2H^2+4(n-1)c}-(n-2)|H|\right).
\]

If $\sup_\Sigma|\Phi|=+\infty$, then (ii) holds trivially and there is nothing to prove. If
$\sup_\Sigma|\Phi|<+\infty$, then by applying \rf{OYw} to the function $|\Phi|^2$ we know that there exists
$\{p_k\}_{k\in\mathbb{N}}$ in $\Sigma$ such that
\[
\lim_{k \rightarrow \infty}|\Phi|(p_k)=\sup_\Sigma|\Phi|,\quad \mathrm{and} \quad \Delta|\Phi|^2(p_k)<1/k,
\]
which jointly with \rf{EQ7} implies
\[
1/k>\Delta|\Phi|^2(p_k)\geq -2|\Phi|^2(p_k)P_H(|\Phi|(p_k)).
\]
Taking limits here, we get $0\geq -2(\sup_\Sigma|\Phi|)^2P_H(\sup_\Sigma|\Phi|)$, that is
\beq
\label{eq2.3}
(\sup_\Sigma|\Phi|)^2P_H(\sup_\Sigma|\Phi|)\geq 0.
\eeq
It follows from here that either $\sup_\Sigma|\Phi|=0$, which means that $|\Phi|=\mathrm{constant}=0$ and the hypersurface
is totally umbilical, or $\sup_\Sigma|\Phi|>0$ and then $P_H(\sup_\Sigma|\Phi|)\geq 0$. In the latter, it must
be $\sup_\Sigma|\Phi|\geq\alpha_H$, which gives the inequality in (ii). Moreover, assume that equality holds,
$\sup_\Sigma|\Phi|=\alpha_H$. In that case, $P_H(|\Phi|)\leq 0$ on
$\Sigma$, which jointly with \rf{EQ7} implies that $|\Phi|^2$ is a subharmonic function on $\Sigma$. Therefore, if
there exists a point $p_0\in\Sigma$ at which this supremum is attained, then $|\Phi|^2$ is a subharmonic
function on $\Sigma$ which attains its supremum at some point of $\Sigma$ and, by the maximum principle, it must
be constant, $|\Phi|=\mathrm{constant}=\alpha_H$. Thus, \rf{EQ7} becomes trivially an equality,
\[
\frac{1}{2}\Delta|\Phi|^2=0=-|\Phi|^2P_H(|\Phi|).
\]
From here we obtain that $\nabla\Phi=\nabla A=0$, that is, the second fundamental form of the hypersurface is
parallel. If $H=0$ (which can occur only when $c=1$) then by a classical local rigidity result by Lawson \cite[Proposition 1]{La}
we know that $\Sigma^n$ is an open piece of a minimal Clifford torus of the form
$\mathbb{S}^k(\sqrt{k/n})\times\mathbb{S}^{n-k}(\sqrt{(n-k)/n})\subset\mathbb{S}^{n+1}$, with $k=1,\ldots,n-1$, which
trivially satisfies
$|\Phi|=\mathrm{constant}=\alpha_0=\sqrt{n}$. If $H\neq 0$ then from the equality in \rf{EQ7} we also obtain
the equality in Okumura lemma (Lemma\rl{okumura}), which implies that the hypersurface
has exactly two constant principal curvatures, with multiplicities $(n-1)$ and $1$.
Then, by the classical results on isoparametric hypersurfaces of Riemannian space forms
\cite{LC,Se,Ca} we conclude that \m\ must be an open piece of one of the three following standard
product embeddings:
\begin{enumerate}
\item[(a)] $\mathbb{R}^{n-1}\times\mathbb{S}^{1}(r)\subset\mathbb{R}^{n+1}$ or
$\mathbb{R}\times\mathbb{S}^{n-1}(r)\subset\mathbb{R}^{n+1}$ with $r>0$, if $c=0$;
\item[(b)] $\mathbb{S}^{1}(\sqrt{1-r^2})\times\mathbb{S}^{n-1}(r)\subset\mathbb{S}^{n+1}$,
with $0<r<1$, if $c=1$; and
\item[(c)] $\mathbb{H}^{n-1}(-\sqrt{1+r^2})\times\mathbb{S}^{1}(r)\subset\mathbb{H}^{n+1}$, with
$0<r<1/\sqrt{n(n-2)}$ (recall that $H^2>-c=1$), or
$\mathbb{H}^{1}(-\sqrt{1+r^2})\times\mathbb{S}^{n-1}(r)\subset\mathbb{H}^{n+1}$, with $r>0$, if $c=-1$.
\end{enumerate}
Obviously, in all the examples above $|\Phi|=\mathrm{constant}=\sup_\Sigma|\Phi|$. A detailed analysis of
the value of the constant $|\Phi|$ for these examples shows that when $c=0$
$|\Phi|=\sqrt{n(n-1)}|H|>\alpha_H$ for the standard products $\mathbb{R}^{n-1}\times\mathbb{S}^{1}(r)$, whereas
$|\Phi|=\sqrt{n}|H|/\sqrt{n-1}=\alpha_H$
for the standard products $\mathbb{R}\times\mathbb{S}^{n-1}(r)$, with $r>0$. On the other hand, when $c=1$ we can see that
\[
|\Phi|=\frac{\sqrt{n}}{2\sqrt{n-1}}\left(\sqrt{n^2H^2+4(n-1)}+(n-2)|H|\right)>\alpha_H
\]
for the standard products $\mathbb{S}^{1}(\sqrt{1-r^2})\times\mathbb{S}^{n-1}(r)$ if $r>\sqrt{(n-1)/n}$, whereas
\[
|\Phi|=\frac{\sqrt{n}}{2\sqrt{n-1}}\left(\sqrt{n^2H^2+4(n-1)}-(n-2)|H|\right)=\alpha_H
\]
if $0<r<\sqrt{(n-1)/n}$. Finally, when $c=-1$ we have that
\[
|\Phi|=\frac{\sqrt{n}}{2\sqrt{n-1}}\left(\sqrt{n^2H^2-4(n-1)}+(n-2)|H|\right)>\alpha_H
\]
for the standard products $\mathbb{H}^{n-1}(-\sqrt{1+r^2})\times\mathbb{S}^{1}(r)$, with $0<r<1/\sqrt{n(n-2)}$,
whereas
\[
|\Phi|=\frac{\sqrt{n}}{2\sqrt{n-1}}\left(\sqrt{n^2H^2-4(n-1)}-(n-2)|H|\right)=\alpha_H
\]
for the standard products $\mathbb{H}^{1}(-\sqrt{1+r^2})\times\mathbb{S}^{n-1}(r)$, with $r>0$.
For the details, see Appendix. This finishes the proof of Theorem\rl{th3a}.

\subsection{Proof of Corollary\rl{coros}}
As in the previous proof, instead of proving Corollary\rl{coros} we will prove its equivalent statement in terms of
the total umbilicity tensor. That is, we will show that, under the assumptions of Corollary\rl{coros}, it holds that
\begin{itemize}
\item[(i)] either $\sup_\Sigma|\Phi|=0$ and \m\ is a totally umbilical hypersurface,
\item[(ii)] or
\[
\sup_\Sigma|\Phi|\geq\alpha_H=\frac{\sqrt{n}}{2\sqrt{n-1}}\left(\sqrt{n^2H^2+4(n-1)c}-(n-2)|H|\right)>0,
\]
and the equality holds and this supremum is attained at some point of $\Sigma$ if and only if one has (a), (b) or (c).
\end{itemize}
Obviously, if $\sup_\Sigma|\Phi|=+\infty$, then (ii) holds trivially and there is nothing to prove. If
$\sup_\Sigma|\Phi|<+\infty$, then we can estimate
\[
H\g{\Phi X}{X}\geq -|H||\g{\Phi X}{X}|\geq -|H||\Phi||X|^2|\geq -|H|\sup_\Sigma|\Phi||X|^2,
\]
and
\[
\g{\Phi X}{\Phi X}\leq |\Phi|^2|X|^2\leq (\sup_\Sigma|\Phi|)^2|X|^2,
\]
for $X\in\mathcal{X}(\Sigma)$. Then, by \rf{Ricci} we obtain for every $X\in \mathcal{X}(\Sigma)$,
\begin{eqnarray*}
\mathrm{Ric}(X,X) & = & (n-1)(c+H^2)|X|^2+(n-2)H\g{\Phi X}{X}-\g{\Phi X}{\Phi X}\\
{} & \geq & \left((n-1)(c+H^2)-(n-2)|H|\sup_\Sigma|\Phi|-(\sup_\Sigma|\Phi|)^2\right)|X|^2.
\end{eqnarray*}
Therefore, if $\sup_\Sigma|\Phi|<+\infty$ then the Ricci curvature of $\Sigma$ is bounded from below by the constant
\[
C=(n-1)(c+H^2)-(n-2)|H|\sup_\Sigma|\Phi|-(\sup_\Sigma|\Phi|)^2.
\]
Since $\Sigma$ is complete, the classical Omori-Yau maximum principle holds on $\Sigma$ and the result follows directly from
Theorem\rl{th3a} (or, equivalently, Theorem\rl{th3aBIS}).

\begin{remark}
The proof of Corollary\rl{coros} has been inspired by the estimates of the Ricci curvature for submanifolds into a Riemannian
space form given by Asperti and Costa in \cite{AsdC}. We refer the reader to that paper for other general Ricci estimates.
\end{remark}

\subsection{Proof of Corollary\rl{coro3a}}
For the proof of Corollary\rl{coro3a}, first recall that the weak Omori-Yau maximum principle holds on every
parabolic Riemannian manifold. Then, by the first part of Theorem\rl{th3a} we obtain that
either $\sup_\Sigma|\Phi|=0$ and \m\ is a totally umbilical hypersurface, or
$\sup_\Sigma|\Phi|\geq\alpha_H$. Moreover, if equality holds, $\sup_\Sigma|\Phi|=\alpha_H$, then as in the proof above
we have $P_H(|\Phi|)\leq 0$ and $|\Phi|^2$ is a subharmonic function on $\Sigma$ which is bounded from above. Since \m\
is parabolic, it must be constant, $|\Phi|=\mathrm{constant}=\alpha_H$. The proof then finishes as in Theorem\rl{th3a},
by observing that the standard Riemannian products $\mathbb{R}\times\mathbb{S}^{n-1}(r)$,
$\mathbb{S}^k(\sqrt{k/n})\times\mathbb{S}^{n-k}(\sqrt{(n-k)/n})$,
$\mathbb{S}^{1}(\sqrt{1-r^2})\times\mathbb{S}^{n-1}(r)$ and
$\mathbb{H}^{1}(-\sqrt{1+r^2})\times\mathbb{S}^{n-1}(r)$ are all parabolic. For
$\mathbb{S}^k(\sqrt{k/n})\times\mathbb{S}^{n-k}(\sqrt{(n-k)/n})$ and
$\mathbb{S}^{1}(\sqrt{1-r^2})\times\mathbb{S}^{n-1}(r)$ this is clear because they are compact.
For the other cases, it follows from the fact that any standard Riemannian product $\mathbb{R}\times M$ with
$M$ compact is parabolic (see \cite[Subsection 2.1]{Ka}).

\section*{Appendix}

In this section we will briefly compute the value of $|\Phi|$ for the standard examples which appears in Theorem\rl{th3a} and
Corollary\rl{coro3a}. In the Euclidean space ($c=0$), apart from the totally umbilical hypersurfaces, the easiest constant mean curvature
hypersurfaces are the standard product embeddings of the form
$\mathbb{R}^{n-k}\times\mathbb{S}^k(r)\hookrightarrow\mathbb{R}^{n+1}$, for a given radius $r>0$ and integer
$k\in\{ 1,\ldots, n-1\}$. Its principal curvatures are given by
\[
\kappa_1=\cdots=\kappa_{n-k}=0, \quad
\kappa_{n-k+1}=\cdots=\kappa_n=\frac{1}{r},
\]
and its constant mean curvature $H$ is given by $nH=k/r$. For these examples
\[
|\Phi|=\frac{\sqrt{k(n-k)}}{\sqrt{n}r}=\frac{\sqrt{n(n-k)}}{\sqrt{k}}|H| \quad
\mathrm{and}
\quad
\alpha_H=\frac{\sqrt{n}}{\sqrt{n-1}}|H|.
\]
In particular, $|\Phi|=\alpha_H$ if and only $k=n-1$, and $|\Phi|>\alpha_H$ otherwise.

When $c=1$, let us consider the standard immersions $\mathbb{S}^{1}(\sqrt{1-r^2})\hookrightarrow\mathbb{R}^{2}$ and
$\mathbb{S}^{n-1}(r)\hookrightarrow\mathbb{R}^{n}$, for a given
radius $0<r<1$, and take the product immersion
$\mathbb{S}^{1}(\sqrt{1-r^2})\times\mathbb{S}^{n-1}(r)\hookrightarrow\mathbb{S}^{n+1}\subset\mathbb{R}^{n+2}$.
Its principal curvatures are given by
\[
\kappa_{1}=\frac{r}{\sqrt{1-r^2}}, \quad \kappa_{2}=\cdots=\kappa_{n}=-\frac{\sqrt{1-r^2}}{r},
\]
and its constant mean curvature is given by
\beq
\label{eqHr}
H=H(r)=\frac{nr^2-(n-1)}{nr\sqrt{1-r^2}}.
\eeq
In this case,
\[
|\Phi|=\frac{\sqrt{n-1}}{r\sqrt{n(1-r^2)}}
\]
where, by \rf{eqHr},
\[
r^2=\frac{2(n-1)+nH^2\pm|H|\sqrt{n^2H^2+4(n-1)}}{2n(1+H^2)}
\]
where we choose the sign $-$ or $+$ according to $r^2\leq (n-1)/n$ or $r^2>(n-1)/n$. Therefore,
\[
|\Phi|=\frac{\sqrt{n}}{2\sqrt{n-1}}\left((n-2)|H|\pm\sqrt{n^2H^2+4(n-1)}\right),
\]
where we use the same criterion for the sign. In particular, $|\Phi|=\alpha_H$ when $r^2\leq (n-1)/n$, and
$|\Phi|>\alpha_H$ when  $r^2>(n-1)/n$.

Finally, when $c=-1$ let us consider the
standard immersions $\mathbb{H}^{n-k}(-\sqrt{1+r^2})\hookrightarrow\mathbb{R}^{n-k+1}_1$ and
$\mathbb{S}^{k}(r)\hookrightarrow\mathbb{R}^{k+1}$, for a given
radius $r>0$ and integer $k\in\{ 1,\ldots, n-1\}$, and take the product immersion
$\mathbb{H}^{n-k}(-\sqrt{1+r^2})\times\mathbb{S}^{k}(r)\hookrightarrow\mathbb{H}^{n+1}\subset\mathbb{R}^{n+2}_1$.
Its principal curvatures are given by
\[
\kappa_{1}=\cdots=\kappa_{n-k}=\frac{r}{\sqrt{1+r^2}}, \quad
\kappa_{n-k+1}=\cdots=\kappa_{n}=\frac{\sqrt{1+r^2}}{r},
\]
and its constant mean curvature is given by
\beq
\label{eqHr2}
H=\frac{nr^2+k}{nr\sqrt{1+r^2}}.
\eeq
We are interested in the cases where $k=1$ and $k=n-1$. Observe that when $k=1$, $H^2>1$ if and only if
$r<1/\sqrt{n(n-2)}$. In that case
\beq
\label{eee}
|\Phi|=\frac{\sqrt{n-1}}{r\sqrt{n(1+r^2)}}
\eeq
where, by \rf{eqHr2},
\[
r^2=\frac{2-nH^2+|H|\sqrt{n^2H^2-4(n-1)}}{2n(H^2-1)}.
\]
Thus,
\[
|\Phi|=\frac{\sqrt{n}}{2\sqrt{n-1}}\left((n-2)|H|+\sqrt{n^2H^2-4(n-1)}\right)>\alpha_H.
\]
On the other hand, when $k=n-1$ we have that $H^2>1$ for every $r>0$ and $|\Phi|$ is also given by \rf{eee}, where now,
by \rf{eqHr2}, $r^2$ is given by
\[
r^2=\frac{2(n-1)-nH^2+|H|\sqrt{n^2H^2-4(n-1)}}{2n(H^2-1)}.
\]
Therefore, in this case we have for every $r>0$
\[
|\Phi|=\frac{\sqrt{n}}{2\sqrt{n-1}}\left((n-2)|H|-\sqrt{n^2H^2-4(n-1)}\right)=\alpha_H.
\]

\bibliographystyle{amsplain}

\end{document}